\documentclass[graybox]{svmult}

\usepackage{type1cm}        
\usepackage{amsmath,amssymb,amsthm}
\usepackage{dsfont,bbm,bm}
\usepackage{makeidx}         
\usepackage{graphicx}        
\usepackage{wrapfig}         
\usepackage{multicol}        
\usepackage{multirow}
\usepackage[bottom]{footmisc}
\usepackage{float}
\usepackage[caption=false]{subfig} 

\newcommand{\STAB}[1]{\begin{tabular}{@{}c@{}}#1\end{tabular}}

\usepackage{newtxtext}       %
\usepackage[varvw]{newtxmath}       

\usepackage{ulem}

\newcommand{\cT}{{\mathcal{T}}}

\makeindex             

\begin{document}

\title*{Adaptive and frugal BDDC coarse spaces for virtual element discretizations of a Stokes problem with heterogeneous viscosity}

\titlerunning{BDDC coarse spaces for virtual element discretizations of Stokes equations}
\author{Tommaso Bevilacqua \and Axel Klawonn \and Martin Lanser}
\institute{Tommaso Bevilacqua$^{1}$, Axel Klawonn$^{1,2}$, Martin Lanser$^{1,2}$ \at $^1$Department of Mathematics and Computer Science, Division of Mathematics, University of Cologne, Weyertal 86-90, 50931 Cologne, Germany, \email{tommaso.bevilacqua@uni-koeln.de, axel.klawonn@uni-koeln.de, martin.lanser@uni-koeln.de}, {url:~\url{https://www.numerik.uni-koeln.de}} \at $^2$Center for Data and Simulation Science, University of Cologne, Germany, {url:~\url{https://www.cds.uni-koeln.de}}}

%
%
\maketitle

\abstract{The virtual element method (VEM) is a family of numerical methods to discretize partial differential equations on general polygonal or polyhedral computational grids. However, the resulting linear systems are often ill-conditioned and robust preconditioning techniques are necessary for an iterative solution. Here, a balancing domain decomposition by constraints (BDDC) preconditioner is considered. Techniques to enrich the coarse space of BDDC applied to a Stokes problem with heterogeneous viscosity are proposed. In this framework a comparison between two adaptive techniques and a computationally cheaper heuristic approach is carried out. Numerical results computed on a physically realistic model show that the latter approach in combination with the deluxe scaling is a promising alternative.}

\section{Introduction}
\label{sec:1}

The virtual element method (VEM) \cite{vem2022book,basic2013} is a finite element discretization approach for partial differential equations (PDEs) which can deal with very general polygonal/polyhedral computational grids. The linear systems that arise from  these discretizations of PDEs are then generally worse conditioned than in the case of  standard finite element methods (FEMs) for which recent studies proposed robust BDDC methods \cite{bertoluzza2017bddc,bertoluzza2020,bevilacqua2022bddc,klawonn2022bddc}.
In this work we analyze a Stokes problem with high heterogeneity in the viscosity function;  thus an adequatly enriched coarse space is needed.
In the adaptive framework  considered here, this is done by solving a generalized eigenvalue problem on each subdomain edge and by adding the solutions to the coarse space in an appropriate way.  In our numerical simulations we make use of two coarse spaces analyzed for standard low order FEM in the paper \cite{klawonn2016comparison}, identifying with "First" the approach present in Section 4.5 and with "Second" the one in Section 5.
The first approach, was originally introduced in \cite{pechsteinmodern} and already successfully used for the  VEM in \cite{bevilacqua2023bddc,dassiZS2022}. The second one, extensively used in dual-primal finite element tearing and interconnecting (FETI-DP) and BDDC in \cite{klawonn2016newton,mandel2007adaptive,mandel2012adaptive,sousedik2013adaptive}, has been also recently extended to the  VEM for diffusion and linear elasticity problems in \cite{klawonn2022adaptive}.\\
An alternative approach to enrich the coarse space denoted as {\it frugal}, has been introduced in \cite{heinlein2019Frugal} and already succesfully used for the VEM for stationary diffusion and linear elasticity in \cite{klawonn2022adaptive}. This is a heuristic and cheaper technique, that does not involve the solution of eigenvalue problems. This often allows to construct robust coarse spaces in a computationally efficient way when it is sufficient to approximate the largest, or smallest, eigenvalues depending on the chosen coarse space.\\
In the present work we extend  the first and second adaptive coarse space approaches as well as the frugal coarse space to the virtual element discretization of a Stokes problem with a heterogeneous viscosity function.

\section{Continuous problem and virtual element discretization}
\label{sec:2}
Let $\mathit{\Omega} \subseteq \mathbb{R}^2$ be a bounded Lipschitz domain, with $\Gamma = \partial \mathit{\Omega}$, and consider the stationary Stokes problem with homogeneous Dirichlet boundary conditions: find $(\mathbf{u},p)$ s.t.  
\begin{equation}\label{ContinuousProblem}
	\begin{cases}
		-\nu(\mathbf{x})\,\mathbf\Delta\mathbf{u} - \nabla p = \mathbf{f} \qquad &\text{ in } \mathit{\Omega} \\
		\text{div}\,\mathbf{u} = 0  \qquad &\text{ in } \mathit{\Omega} \\
		\mathbf{u}=0 \qquad  &\text{ on } \mathit{\Gamma}.
	\end{cases}
\end{equation}
Here $\mathbf{u}$ and $p$ are respectively the velocity and the pressure fields, $\mathbf{f}\in [H^{-1}(\mathit{\Omega})]^2$ represents the external force and $\nu\in L^{\infty}(\mathit{\Omega})$, $\nu(\mathbf{x}) >0$, $\mathbf{x} \in \mathit{\Omega}$, is the heterogeneous viscosity function. \\
Introducing $\mathbf{V}:=[H^{1}_{0}(\mathit{\Omega})]^2$ and $Q:=L^{2}_{0}(\mathit{\Omega})= \{ q\in L^{2}(\mathit{\Omega})\, s.t. \, \int_{\mathit{\Omega}} q =0 \}$, the standard variational formulation reads: find $(\mathbf{u},p) \in \mathbf{V} \times Q$ s.t.
\begin{equation}\label{VarForm}
	\begin{cases}
		a(\mathbf{u},\mathbf{v})+b(\mathbf{v},p)=(\mathbf{f},\mathbf{v}) & \text{for all } \mathbf{v} \in \mathbf{V}, \\
		b(\mathbf{u},q)=0 & \text{for all } q \in Q,
	\end{cases}
\end{equation}
where $\,a(\mathbf{v},\mathbf{w}) := \int_{\mathit{\Omega}}\nu(\mathbf{x}) \mathbf{\nabla v} : \mathbf{\nabla w}$ for all $\mathbf{v},\mathbf{w} \in \mathbf{V}$, $\,b(\mathbf{v},q) := \int_{\mathit{\Omega}} \text{div}\, \mathbf{v} \, q $ for all $\mathbf{v} \in \mathbf{V}$, $q \in Q$ and $\,(\mathbf{f},\mathbf{v}):=\int_{\mathit{\Omega}} \mathbf{f} \cdot \mathbf{v}$  for all $\mathbf{v} \in \mathbf{V}$.

The discretization of problem \eqref{VarForm} is based on a virtual element space which is designed to solve a Stokes problem. In the following we present the basic elements of this discretization; we refer  to \cite{da2017divergence} for further details.\\
Let $\{\cT_h\}_h$ be a sequence of triangulations of $\mathit{\Omega}$ into general polygonal elements $K$ with $	h_K:=\text{diameter}(K)$ and $h:=\sup_{K\in \cT_h}h_K$.
We suppose that, for all $h$, each element $K\in \cT_h$ satisfies, for some $\gamma>0$ and $c>0$, the following assumptions
\begin{itemize}
\itemsep0em 
	\item $K$ is star-shaped with respect to a ball of radius greater or equal than $\gamma h_k$,
	\item  the distance between any two vertices of $K$ is greater or equal than $c h_K$,
\end{itemize} 

For $k \in \mathbb{N}$, we then define the spaces: $\mathbb{P}_k(K)$ the set of polynomials on $K$ of degree smaller or equal than $k$, $\mathbb{B}_k(K):=\{v \in C^0(\partial K ) \text{ s.t. } v_{|e} \in \mathbb{P}_k(e) \quad \forall \text{ edge } e \in \partial K \}$, $\mathit{G}_k(K):=\nabla(\mathbb{P}_{k+1}(K)) \subseteq [\mathbb{P}_{k}(K)]^2$ and its $L^2$-orthogonal complement $\mathit{G}_k(K)^\perp \subseteq [\mathbb{P}_{k}(K)]^2$.
The local virtual element spaces are defined, for $k\geq2$, on each $K \in \mathit{T}_h$ as
\begin{eqnarray*}\label{NewLoc}
	& \mathbf{V}_h^K:=\bigg\{ \mathbf{v} \! \in \! [H^1(K)]^2 \,\big |\, \mathbf{v}_{| \partial K} \in [\mathbb{B}_{k}(\partial K)]^2,
		\begin{cases} -\nu\mathbf\Delta\mathbf{v} - \nabla s \in\mathit{G}_{k-2}(K)^\perp, \\ \text{div } \mathbf{v} \in \mathbb{P}_{k-1}(K),\end{cases} \mkern-20mu s \in L^2(K)\bigg\},\\
		 & Q_h^K:=\mathbb{P}_{k-1}(K),
\end{eqnarray*}
and the global virtual element spaces are	$\mathbf{V}_h:=\{\mathbf{v}\in [H^{1}_{0}(\mathit{\Omega})]^2 \, \big | \, \mathbf{v}_{|K} \in \mathbf{V}^K_h, \, \forall K\in \cT_h\}$ and $Q_h:=\{q\in L^{2}_{0}(\mathit{\Omega})  \, \big | \, q_{|K} \in Q^K_h, \, \forall K\in \cT_h\}$.\\
In the VEM framework the basis function are never explicitly computed since it would be necessary solve the PDE given in $\mathbf{V}_h^K$ for each element. Alternatively, they are defined by using of some polynomial projection operators (see \cite{ahmad2013equivalent}) and suitable degrees of freedom (dofs).

Given $\mathbf{v} \in \mathbf{V}_h^K$ we take the following linear operators $\mathbf{D_{V}}$, split into 
\begin{itemize}
\itemsep0em 
	\item $\mathbf{D_{V}^1}$: the values of $\mathbf{v}$ at the vertices of the polygon $K$,
	\item $\mathbf{D_{V}^2}$: the values of $\mathbf{v}$ at $k-1$ internal points of the $(k+1)$-Gauss-Lobatto quadrature rule in $e\in \partial K$,
	\item $\mathbf{D_{V}^3}$: the moments of $\mathbf{v}$: $\int_{K}\mathbf{v} \cdot \mathbf{g}^\perp_{k-2} \quad \text{for all } \mathbf{g}^\perp_{k-2} \in \mathbf{G}_{k-2}(K)^\perp$;
	\item $\mathbf{D_{V}^4}$: the moments of div$\, \mathbf{v}$: $\int_{K} \text{div} \, \mathbf{v} \,q_{k-1} \quad 	\text{for all } q_{k-1}\in \mathbb{P}_{k-1}(K)/\mathbb{R}$.
\end{itemize} 
Furthermore, for the local pressure, given $q\in Q^K_h$, 
 the linear operators $\mathbf{D_Q}$:
\begin{itemize}
	\item $\mathbf{D_Q}$: the moments of $q$:
	$\int_{K}q p_{k-1} ,	\quad \text{for any } p_{k-1} \in \mathbb{P}_{k-1}(K)$.
\end{itemize}
We note that $b(\mathbf{v}_h,q_h)$ for all $\mathbf{v}_h \in \mathbf{V}_h$ and $q_h \in Q_h$ can be computed directly from the $\mathbf{D^1_V}, \mathbf{D^2_V}$ and $\mathbf{D^3_V}$. While $a(\mathbf{v}_h,\mathbf{w}_h)$ and $(\mathbf{f},\mathbf{v}_h)$ for all $\mathbf{v}_h, \mathbf{w}_h \in \mathbf{V}_h$ can
not be exactly computed. It is therefore necessary to introduce approximations $a_h(\mathbf{v}_h,\mathbf{w}_h)$ and $(\mathbf{f}_h,\mathbf{v}_h)$ making use of suitable polynomial projections. Further details of the construction of these bilinear forms and their related theoretical estimates can be found in \cite{da2017divergence}.
The discrete virtual element problem states: find $(\mathbf{u}_h,p_h)$ s.t.
\begin{equation}\label{DiscreteProblem}
	\begin{cases}
		a_h(\mathbf{u}_h,\mathbf{v}_h)+b(\mathbf{v}_h,p_h)=(\mathbf{f}_h,\mathbf{v}_h)& \quad \text{for all }\mathbf{v}_h \in \mathbf{V}_h\\
		b(\mathbf{u}_h,q_h) = 0&\quad \text{for all }q_h \in Q_h
	\end{cases}.
\end{equation}

\section{Domain decomposition, BDDC preconditioner, and coarse spaces}
We decompose $\mathcal{T}_h$ into $N$ non-overlapping subdomains $\mathit{\Omega}_i$ with characteristic size $H_i$ as $\bar{\mathcal{T}}_h = \bigcup_{i=1}^N \bar{\mathit{\Omega}}_i$ where each $\mathit{\Omega}_i$ is the union of different polygons of the tessellation  $\mathcal{T}_h$ and we define $\Gamma = \bigcup_{i\neq j} \partial \mathit{\Omega}_i \cap \partial \mathit{\Omega}_j$ as interface among the subdomains. We assume that the decomposition is shape-regular in the sense of \cite{bertoluzza2017bddc} Section 3.\\
We refer to the edges of the subdomains $\mathit{\Omega}_i$ as macro edges and we denote them with $\mathcal{E}_i$, moreover $\mathcal{E}_{ij}$ denote the macro edge shared by the subdomains $\mathit{\Omega}^i$ and $\mathit{\Omega}^j$.

From now we omit the underscore $h$ since we will always refer to the finite-dimensional space and we write $\mathbf{V}\times Q$ instead of $\mathbf{V}_{h}\times Q_{h}$.
We split the velocity dofs into interface ($\Gamma$) and internal ($I$) dofs. In particular, the $\mathbf{D}^1_V$ and $\mathbf{D}^2_V$ that belongs to a single subdomain $\mathit{\Omega}_i$ are classified as internal, while the ones that belong to more than a single subdomain as interface dofs. All the dofs $\mathbf{D}^3_V$ and $\mathbf{D}^4_V$ are classified as internal ones.
Following the notations introduced in \cite{li2006bddc} and \cite{bevilacqua2022bddc}, we decompose the discrete velocity and pressure space $\mathbf{V}$ and $Q$ into $	\mathbf{V} = \mathbf{V}_I \bigoplus \mathbf{\widehat{V}}_\Gamma$ and $Q = Q_I\bigoplus Q_0$, with $Q_0:=\prod_{i=1}^N \{q\in \mathit{\Omega}_i | q \textit{ is constant in } \mathit{\Omega}_i\}$.\\
$\mathbf{\widehat{V}}_\Gamma$ is the continuous space of the traces on $\Gamma$ of functions in $\mathbf{V}$, $\mathbf{V}_I = \bigoplus_{i=1}^{N} \mathbf{V}_I^{(i)}$ and $Q_I = \bigoplus_{i=1}^{N} Q_I^{(i)}$ are direct sums of subdomain interior velocity and pressure spaces.

We also define the space of interface velocity variables of the subdomain $\mathit{\Omega}_i$ by $\mathbf{V}_\Gamma^{(i)}$ and the associated product space by $\mathbf{V}_\Gamma = \prod_{i=1}^{N}\mathbf{V}_\Gamma^{(i)}$.
The discrete global saddle-point problem \eqref{DiscreteProblem} can be written as: find $(\mathbf{u}_I,p_I,\mathbf{u}_\Gamma,p_0) \in (\mathbf{V}_I,Q_I,\mathbf{\widehat{V}}_\Gamma,Q_0)$ s.t.
\begin{equation}\label{discMat}
	\left[
	\begin{array}{cccc}
		A_{II} & B_{II}^T & \widehat{A}_{\Gamma I}^T & 0\\
		B_{II} & 0 & \widehat{B}_{I\Gamma} & 0\\
		\widehat{A}_{\Gamma I} & \widehat{B}_{I\Gamma}^T & \widehat{A}_{\Gamma \Gamma} & \widehat{B}_{0\Gamma}^T\\
		0 & 0 & \widehat{B}_{0\Gamma}^T & 0\\
	\end{array}
	\right]
	\left[
	\begin{array}{c}
		\mathbf{u}_I\\
		p_I\\
		\mathbf{u}_\Gamma\\
		p_0\\
	\end{array}
	\right]
	=
	\left[
	\begin{array}{c}
		\mathbf{f}_I\\
		0\\
		\mathbf{f}_\Gamma\\
		0\\
	\end{array}
	\right],
\end{equation}
where the blocks $\widehat{\cdot}$ related to the continuous interface velocity are assembled from the corresponding subdomain submatrices. By static condensation one eliminates the interior variables and obtains the global interface saddle point problem
\begin{equation}\label{globInt}
	\widehat{S}\text{ }\widehat{u} = 
	\left[
	\begin{array}{cc}
		\widehat{S}_\Gamma & {\widehat{B}_{0\Gamma}}^T\\
		\widehat{B}_{0\Gamma} & 0\\
	\end{array}
	\right]
	\left[
	\begin{array}{c}
		\mathbf{u}_\Gamma\\
		p_0\\
	\end{array}
	\right]
	=
	\left[
	\begin{array}{c}
		\mathbf{g}_\Gamma\\
		0\\
	\end{array}
	\right]
	= \widehat{\mathbf{g}}\text{,}
\end{equation}
where $\widehat{S}_\Gamma$ is the Schur complement of the submatrices constituted by the top $3 \times 3$ block of the left-hand side matrix in \eqref{discMat} and $\widehat{\mathbf{g}}$ the corrispective right-hand side.

We introduce $\mathbf{\widetilde{V}}_\Gamma = \mathbf{\widehat{V}}_\Pi \bigoplus \mathbf{V}_\Delta = \mathbf{\widehat{V}}_\Pi \bigoplus \big( \prod_{i=1}^N \mathbf{V}_\Delta^{(i)} \big)$, a partially assembled interface velocity space where $\mathbf{\widehat{V}}_\Pi$ is the continuous coarse-level primal velocity space, which dofs are shared by neighboring subdomains and the complementary space $\mathbf{V}_\Delta$, that is the direct sum of the subdomain dual interface velocity spaces $\mathbf{V}_\Delta^{(i)}$.\\
In particular, the primal dofs of our problem are represented by the nodal evaluation of both components of the velocity in the vertices of the subdomain and one extra dof for each subdomain edge to satisfy the \textit{no-net-flux} condition $\int_{\partial\mathit{\Omega}_i}\mathbf{v}_\Delta^{(i)}\cdot\mathbf{n} = 0 \,\, \forall \mathbf{v}_\Delta\in\mathbf{V}_\Delta$, where $\mathbf{n}$ is the outward normal of $\partial\mathit{\Omega}_i$. This last condition is needed to ensure that the operator of the preconditioned system \eqref{globInt} with the BDDC is symmetric and positive definite on some particular subspaces \cite{li2006bddc}, so the preconditioned conjugate gradient (CG) method can be used for solution.\\
In our study we use a generalized transformation of basis approach \cite{klawonn2020coarse} such that each primal basis function corresponds to an explicit dof. Firstly, on each $\mathcal{E}_{ij}$, we assume that the velocity vector should fulfill $N_{ij}$ constraints $c_{ij}^{l}$ for $l = 1,...,N_{ij}$, i.e., $c_{ij}^{l^T}\mathbf{v}_{\Gamma|\mathcal{E}}^{(i)} = c_{ij}^{l^T}\mathbf{v}_{\Gamma|\mathcal{E}}^{(j)}$. Then, we compute the orthonormal trasformation with a modified Gram-Schmidt algorithm. Finally, since the transformations are independent of each other, we construct the resulting block diagonal global transformation. The constraints $c_{ij}^{l}$ are established via the \textit{no-net-flux} condition and the techniques to enrich the coarse space. More details  can be found in \cite{klawonn2020coarse}.\\
For each $\mathit{\Omega}_i$ we introduce the scaling matrices $D$. These can be chosen in different ways and can be either diagonal or not, but they always must provide a partition of unity, i.e., $\widetilde{R}_{D,\Gamma}^T \widetilde{R}_\Gamma = \widetilde{R}_\Gamma^T \widetilde{R}_{D,\Gamma} = I$, where $\widetilde{R}_\Gamma,\widetilde{R}_{D,\Gamma} :\mathbf{\widehat{V}}_\Gamma\rightarrow\mathbf{\widetilde{V}}_\Gamma$ are respectively a restriction operator, and its scaled version. We use the standard \textit{multiplicity-scaling} and a variant of the \textit{deluxe-scaling} to preserve the normal fluxes \cite{zampini2017multilevel}.\\
We then define the average operator $E_D = \widetilde{R}\widetilde{R}_D^T$, which maps $\mathbf{\widetilde{V}}_\Gamma\times Q_0$, with generally discontinuous interface velocities, to elements with continuous interface velocities in the same space. Here $\widetilde{R}$ and $\widetilde{R}_D^T$ are simply the two previous operators extended by identity to the space of piecewise constant pressures.\\
The preconditioner for solving the global saddle-point problem \eqref{globInt} is then $M^{-1}=\widetilde{R}_D^T \widetilde{S}^{-1} \widetilde{R}_D,$ where $\widetilde{S}$ is the Schur complement system that arises using the partially assembled velocity interface functions. Theoretical estimates show that the condition number is bounded by the norm of the average operator \cite{bevilacqua2022bddc,li2006bddc}.\\

\noindent
{\bf Adaptive and frugal coarse spaces}\\
The coarse spaces is alternatively enriched by two adaptive techniques or a heuristic one. The idea is to detect the largest, or smaller, eigenvalues on each macro edge and then include the corresponing eigenvectors in the coarse space as primal constraints $c_{ij}$ with the transormation of basis approach we saw before. \\
We define here the frugal coarse space for our model problem. Like in the linear elasticity case, when applying the BDDC method to the Stokes problem in two dimensions, we need three constraints for each edge to control the three (linearized) rigid-body motions (two translations and one rotation). Given two subdomains $\mathit{\Omega}_l$, $l=i,j$ with diameter $H_l$, we have
\begin{align}
\mathbf{r}_1:=
\left[
\begin{array}{c}
1 \quad 0\\
\end{array}
\right],
\quad
\mathbf{r}_2:=
\left[
\begin{array}{c}
0 \quad 1\\
\end{array}
\right],
\quad
\mathbf{r}_3:=\frac{1}{H_l}
\left[
\begin{array}{c}
x_2 - \widehat{x}_2 \quad -x_1 + \widehat{x}_1\\
\end{array}
\right],
\end{align}
where $\widehat{\mathbf{x}} \in \mathit{\Omega}_l$ is the center of the rotation. Differently from the approach in \cite{klawonn2022adaptive}, we do not rescale the rigid body modes and we define the "approximate" eigenvector
\begin{align}
    v(\mathbf{x})^{(m,l)}_{\mathcal{E}_{ij}}:=
    \begin{cases}
        r(\mathbf{x})_m^{(l)}, & \mathbf{x}\in \mathcal{E}_{ij},\\
        0 , & \mathbf{x}\in \partial\mathit{\Omega}_l \setminus \mathcal{E}_{ij},
    \end{cases}
\end{align}
for $m=1,2,3$ and $v(\mathbf{x})^{(m)\,T}_{\mathcal{E}_{ij}}:=[v(\mathbf{x})^{(m,i)\,T}_{\mathcal{E}_{ij}}, -v(\mathbf{x})^{(m,j)\,T}_{\mathcal{E}_{ij}}]$. \\
The three frugal edge constraints are then obtained by $c_{ij} := B_{D_{ij}}^T S_{ij} P_{D_{ij}} \, v(\mathbf{x})^{(m)\,T}_{\mathcal{E}_{ij}}$, where $S_{ij} := diag(S_\Gamma^{(i)},S_\Gamma^{(j)})$, $B_{D_{ij}}$ is a scaled jump operator, $B_{D,\mathcal{E}_{ij}}$ its restriction to the edge $\mathcal{E}_{ij}$, and $P_{D_{ij}} = B_{D,\mathcal{E}_{ij}}^T B_{\mathcal{E}_{ij}}$.

\section{Numerical Results}
\label{sec:6}

\begin{table}[b]
\footnotesize
\centering
\caption{Test with nSink = 11 and increasing the number of subdomains. Mesh sizes: 1000 elements for $2\times2$ subdomains, 5000 elements for $4\times4$ and 12000 elements for $8\times8$.}\label{tab:Sinkers}
\vspace{0.5mm}
\begin{tabular}{|c c|c c c|c c c|c c c|c c c|}
\hline

\multicolumn{2}{|c|}{} & \multicolumn{6}{c|}{muliplicity-scaling} & \multicolumn{6}{c|}{deluxe-scaling}\\
\multicolumn{2}{|c|}{} & \multicolumn{3}{c}{CVT} & \multicolumn{3}{c|}{RND} & \multicolumn{3}{c}{CVT} & \multicolumn{3}{c|}{RND}\\

& & $n_\Pi$ & it  & $k_2$ & $n_\Pi$ & it  & $k_2$ & $n_\Pi$ & it & $k_2$ & $n_\Pi$ & it & $k_2$ \\ \hline
 \multirow{3}{*}{\STAB{\rotatebox[origin=c]{90}{\textit{Frugal}}}}
&2x2  & 24 & 133  & 1538.25 & 24 & 102  & 927.72 & 24  & 9 & 16.47  & 24  & 10 & 10.26 \\  
&4x4  & 168 & 62  & 417.12  & 168 & 93 & 3925.17 & 168 & 16 & 18.20 & 168 & 15 & 38.34 \\  
&8x8  & 842 & 52  & 124.72  & 830 & 53 & 204.71  & 842 & 16 & 10.89 & 830 & 14 & 5.52 \\  
\hline
\hline
\multirow{3}{*}{\STAB{\rotatebox[origin=c]{90}{\textit{First}}}}
&2x2   & 49 & 58 & 106.54   & 56 & 48 & 78.09      & 11 & 11  & 9.92  & 12 & 11 & 10.99 \\      
&4x4   & 96 & 52  & 81.48   & 104 & 54  & 90.37  & 81 & 18  & 10.59   & 84 & 19 & 16.49  \\      
&8x8   & 418 & 54 & 102.54  & 424 & 53 & 98.16 & 386 & 21 & 25.95 & 380 & 21 & 16.93  \\      
\hline
\hline
\multirow{3}{*}{\STAB{\rotatebox[origin=c]{90}{\textit{Second}}}}
&2x2  & 47 & 57   & 107.74   & 33 & 50  & 94.27    & 9 & 11 & 10.11 & 9 & 11 & 11.00 \\ 
&4x4  & 84 & 52 & 103.63 	& 89 & 56 & 97.76   & 69 & 17 & 10.11    & 69 & 20  & 17.09 \\
&8x8  & 387 & 57  & 97.97   & 396 & 56 & 96.37  & 361 & 20  & 10.14 & 362 & 21 & 19.83 \\  	
\hline
\end{tabular}
\end{table}

\begin{wrapfigure}{tr}{0.4\textwidth}
\vspace{-1.25cm}
  \begin{center}
    \includegraphics[width=0.33\textwidth]{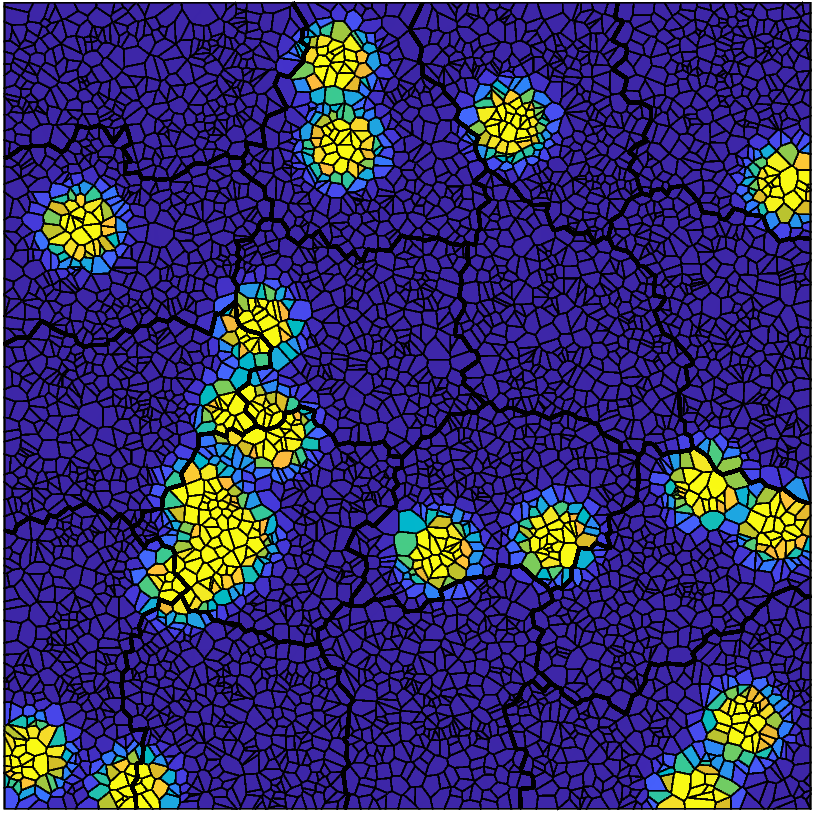}
  \end{center}
\vspace{-0.5cm}
	\caption{METIS decomposition of a RND mesh. $\nu = 1e3$ on yellow elements, $\nu = 1e-3$ on blue ones.}\label{fig:MetisSinkers}
\end{wrapfigure}

We solve a lid-driven cavity benchmark problem on the unit square domain $\mathit{\Omega} = [0,1]\times[0,1]$, applying Dirichlet boundary conditions on the whole $\partial\mathit{\Omega}$ and using a VEM implementation of degree $k = 2$. 
The heterogeneity is introduced to physically represent a practical example where drops (or sinkers) of a high viscosity material are spread in the fluid, in particular this is modeled defining $\nu(\mathbf{x})$ as a continuous function that exhibits sharp gradients (Figure \ref{fig:MetisSinkers}) \cite{rudi2017weighted}. 
These inclusions of equal size are placed randomly in the unit square domain so that they can overlap and intersect the boundary. For $\mathbf{x}\in \mathit{\Omega}$, the viscosity $\nu(\mathbf{x}) \in [\nu_{min},\nu_{max}]$, $0<\nu_{min}<\nu_{max}<\infty$, is defined as $\nu(\mathbf{x}) := (\nu_{max}-\nu_{min})(1-\chi_n(\mathbf{x}))+\nu_{min}$. Here, $\chi_n(\mathbf{x}) \in C^\infty$ is an indicator function $\chi_n(\mathbf{x})\in [0,1]$ that accumulates $n$ sinkers defined as $\chi_n(\mathbf{x}) := \prod_{i=1}^n 1-\exp\big (-\delta \max\big(0,|\mathbf{c}_i-\mathbf{x}|-\frac{\mathit{\Omega}}{2}\big)^2\big )$, where $\mathbf{c}_i \in \mathit{\Omega}, i = 1,...,n$ are the centers of the sinkers, $\mathit{\Omega} \geq 0$ is their diameter and  $\delta>0$ a parameter that controls the exponential decay. By choosing $\delta = 2000$, $\mathit{\Omega} = 0.05$, $\lambda_{min} = 10^{-3}$ and $\lambda_{max} = 10^{3}$, we ensure that the viscosity exibits sharp gradients. 
The right hand side is defined as $\mathbf{f}(\mathbf{x}):=(0,\beta(\chi_n(\mathbf{x}-1)))$, with $\beta = 10$ to simulate gravity that takes down the high viscosity material.
In our experiments we use meshes with a Centroid Voronoi Tassellation (CVT) and Random meshes (RND), while the subdomain partitioning is performed by METIS. We compare the two adaptive coarse spaces, with TOL = 100, and the frugal one by applying the two different type of scaling mentioned before. Our numerical simulation have been performed with MATLAB R2023A\copyright \, therefore no CPU time analysis is provided.

In the following tables, we report the number of iterations to solve the global interface saddle-point problem (\ref{globInt}) with the PCG method, accelerated by a BDDC preconditioner, where we set the tolerance for the relative residual error to $10^{-6}$. \\
In the tables we use the following notation: nSub = number of subdomains, nSink = number of sinkers, $n_\Pi$~=~number of primal constraints, it = iteration count (CG), $k_2$ = condition number, $Frugal$ = frugal coarse space, $First$ =  first adaptive technique, $Second$ = second adaptive technique. 

\begin{table}[t]
\footnotesize
\centering
\caption{Test with increasing number of randomly placed sinkers. Both meshes are made of  4096 elements decomposed into $4\times4$ subdomains.}\label{tab:SinkersInc}
\vspace{0.5mm}
\begin{tabular}{|c c|c c c|c c c|c c c|c c c|}
\hline
\multicolumn{2}{|c|}{} & \multicolumn{6}{c|}{muliplicity-scaling} & \multicolumn{6}{c|}{deluxe-scaling}\\
\multicolumn{2}{|c|}{} & \multicolumn{3}{c}{CVT} & \multicolumn{3}{c|}{RND} & \multicolumn{3}{c}{CVT} & \multicolumn{3}{c|}{RND}\\
& nSink & $n_\Pi$ & it  & $k_2$ & $n_\Pi$ & it  & $k_2$ & $n_\Pi$ & it  & $k_2$ & $n_\Pi$ & it  & $k_2$ \\ \hline
 \multirow{4}{*}{\STAB{\rotatebox[origin=c]{90}{\textit{Frugal}}}}
&1  & 166 & 32 & 60.15     & 168 & 27  & 42.41   & 166 & 10 & 2.54  & 168 & 10  & 6.41  \\
&5  & 166 & 50 & 259.22    & 168 & 85  & 920.27  & 166 & 12 & 35.92 & 168 & 12  & 24.34 \\
&10 & 166 & 81 & 1535.32   & 168 & 118 & 1086.15 & 166 & 14 & 41.66 & 168 & 13  & 50.45 \\
&20 & 166 & 106 & 2537.60  & 168 & 150 & 1523.15 & 166 & 18 & 96.97 & 168 & 13  & 59.60  \\
\hline
\hline
\multirow{4}{*}{\STAB{\rotatebox[origin=c]{90}{\textit{First}}}}
&1   & 70 & 37  & 85.74   & 70 & 28  & 42.42   & 70 & 15    & 4.97  & 70 & 17  & 4.72 \\
&5   & 85 & 46  & 85.74   & 96 & 51  & 84.01   & 75 & 15    & 5.09  & 75 & 15  & 7.68 \\
&10  & 99 & 54  & 96.28   & 117 & 57 & 98.86   & 80 & 16    & 9.56  & 78 & 15  & 8.32 \\
&20  & 121 & 57 & 88.27   & 163 & 66 & 135.31  & 89 & 19    & 13.14 & 81 & 20  & 10.73 \\
\hline
\hline
\multirow{4}{*}{\STAB{\rotatebox[origin=c]{90}{\textit{Second}}}}
&1    & 70 & 34 & 85.68     & 70 & 30    & 42.42   & 71 & 15 &4.82   & 69 & 17      & 4.72  \\
&5    & 80 & 46 & 86.96     & 93 & 52    & 91.34   & 75 & 14 &5.49   & 74 & 15      & 9.97  \\
&10   & 89 & 56 & 97.06     & 107 & 59   & 98.83   & 69 & 17 &9.57   & 69 & 17      & 10.54 \\
&20   & 106 & 57 & 94.27     & 150 & 64 &98.81     & 73 & 20 & 15.79  & 71 & 21    & 10.84  \\
\hline
\end{tabular}
\end{table}

We consider two different tests. We first set a configuration with nSink = 11 and we increase the number of the subdomains; see Table \ref{tab:Sinkers}. For both the type of the mesh considered the adaptive coarse spaces with the multiplicity scaling respect our expectations, while the frugal approach exibits a high condition number since the heuristic coarse space is not able to catch all the largest eigenvalues. Introducing the deluxe scaling we see that the number of primal constraints is drastically reduced in the adaptive coarse spaces. The frugal approach is then able to control the largest eigenvalues and performs well. 
In Table \ref{tab:SinkersInc} we instead keep fixed the number of subdomains at $4\times4$ and we increase the number of the inclusions. Again, the adaptive coarse spaces are robust and also when introducing the deluxe scaling the frugal one shows a good improvement presenting itself as a valid alternative.

\bibliographystyle{spmpsci}
\bibliography{literature}
\end{document}